\documentclass{article}
\usepackage{amsmath,amstext,amsthm,amsfonts,latexsym}
\usepackage[dvips]{graphicx}
\usepackage{epic,eepic}
\theoremstyle{plain}
\newtheorem{theorem}{Theorem }[section]
\newtheorem{proposition}[theorem]{Proposition}
\newtheorem{lemma}[theorem]{Lemma}

\newtheorem{maintheorem}{Theorem}

\theoremstyle{definition}
\newtheorem{remark}[theorem]{Remark}

\newtheorem{definition}[theorem]{Definition}

\newtheorem{question}{Question}

\newcommand{\field}[1]{\mathbb{#1}}
\newcommand{\real}{\field{R}}

\newcommand{\integer}{\field{Z}}

\newcommand{\Tor}{\field{T}}

\newcommand{\al} {\alpha}       
        
\newcommand{\ga} {\gamma}    
\newcommand{\de} {\delta}       \newcommand{\De}{\Delta}
     
\newcommand{\vep}{\varepsilon}
\newcommand{\ze} {\zeta}

\newcommand{\la} {\lambda}      \newcommand{\La}{\Lambda}

\newcommand{\om} {\omega}

\newcommand{\SA}{{\cal A}}
\newcommand{\SB}{{\cal B}}

\newcommand{\p}{\partial}
\newcommand{\supp}{\text{supp}}
\newcommand{\infint}{\int_{-\infty}^{\infty}}

\begin{document}

\title{On the periodic Schr\"odinger-Debye equation}

\author{ Alexander Arbieto and Carlos Matheus\footnote{The authors were
partially supported by Faperj-Brazil.} }

\date{November 24, 2005}

\maketitle

\begin{abstract}We study local and global well-posedness of the initial value 
problem for the Schr\"odinger-Debye equation in the \emph{periodic case}. More
precisely, we prove local well-posedness for the periodic Schr\"odinger-Debye
equation with subcritical nonlinearity in arbitrary dimensions. Moreover, we
derive a new \emph{a priori} estimate for the $H^1$ norm of solutions of the
periodic Schr\"odinger-Debye equation. A novel phenomena obtained as a
by-product of this \emph{a priori} estimate is the global well-posedness of the
periodic Schr\"odinger-Debye equation in dimensions $1,2$ and $3$ \emph{without}
any smallness hypothesis of the $H^1$ norm of the initial data in the 
``focusing'' case.
\end{abstract}

\section{Introduction}

The main theme of this paper is the well-posedness of the initial value problem 
(IVP) for the \emph{Schr\"odinger-Debye equation} (SDE):
\begin{eqnarray}\label{e.S-D}
\left\{ \begin{array}{ll}
i\p_t u+\De u=uv,& t\geq 0,\ \ x\in\Tor^n,\\
K\p_t v+v=\vep|u|^{\al},\\
u(0,x)=u_0(x),& v(0,x)=v_0(x),
\end{array}
\right.
\end{eqnarray}
where $\al=p-2>0$, $u$ is a complex-valued function, $v$ is a real-valued function, $K>0$,
$\vep=\pm 1$ and $\De$ is the Laplacian operator in the $x$-variable. 

This equation appears naturally in certain \emph{nonlinear optics} phenomena.
Indeed, the equation~(\ref{e.S-D}) is obtained from the
\emph{Maxwell-Debye system}
\begin{displaymath}
\left\{ \begin{array}{ll}
i\p_t A+\frac{c}{k\eta_0}\De A=\frac{\om_0}{\eta_0}\nu A,\\
K\p_t\nu+\nu=\eta_2|A|^{\al},
\end{array}
\right.
\end{displaymath}
via the rescaling
\begin{displaymath}
\begin{array}{ll}
u(t,x)=\sqrt{\frac{\om_0|\eta_2|}{\eta_0}}A(t,\sqrt{\frac{c}{k\eta_0}}x),\\
v(t,x)=\frac{\om_0}{\eta_0}\nu(t,\sqrt{\frac{c}{k\eta_0}}x).
\end{array}
\end{displaymath}

Physically, the Maxwell-Debye system (with $\al=2$) arises in nonlinear optics 
describing the non-resonant delayed iteraction of an electromagnetic wave with a 
certain media. In this system, $A$ denotes the envelope of a light wave that
travels through a media. The wave induces a change $\nu$ of the refractive index
in the material (initially $\eta_0$ for an electromagnetic wave of frequency
$\om_0$) with a slight delay $K$. The parameter $\eta_2$ is related with the
magnitude and the sign of the coupling of the wave and the matter. Finally, $c$
is the light velocity in the vacuum and $k$ is the wave vector of the incident
electromagnetic wave. See~\cite{CL} and 
references therein for discussions of this model.

Mathematically, the well-posedness of the IVP~(\ref{e.S-D}) \emph{in the
non-periodic case} (i.e., $x\in\real^n$) was recently studied by
Bid\'egaray~\cite{Bi1},~\cite{Bi2} and 
Corcho, Linares~\cite{CL}.

The results proved by Bid\'egaray, roughly speaking, were local well-posedness
in $L^2(\real^n)$ for data $u_0,v_0\in L^2(\real^n)$ and local well-posedness in
$H^1(\real^n)$ for data $u_0,v_0\in H^1(\real^n)$ ($n=1,2,3$), although the
persistence property was not obtained.  

More recently, Corcho and Linares~\cite{CL}, making an optimal 
use of \emph{Strichartz's inequalities} for the linear Schr\"odinger operator, 
were able to improve Bid\'egaray's results (and, in fact, obtain new ones).

The strategy used by Bid\'egaray and Corcho, Linares was a combination of the
Strichartz inequality and a fixed point argument.

In this paper we apply the same strategy to the IVP~(\ref{e.S-D}), namely, use a
fixed point argument and Strichartz inequality in the periodic case. However,
there is an extra difficulty in our case because the exact analogue of the 
Strichartz inequality does not hold in
the torus $\Tor^n$, and Strichartz-like inequalities may only holds locally in
time.
 
The idea to overcome this is to use the work of Bourgain~\cite{B}, where the
correct analogue of Strichartz's inequality was found and applied to the
nonlinear periodic Schr\"odinger equation (again by a fixed point argument).

In order to apply the fixed point method to solve the Schr\"odinger-Debye
equation, we start by decoupling the equation~(\ref{e.S-D}) to obtain the 
integral formulation:

\begin{equation}\label{e.0}
v(t)=e^{-t/K}v_0(x)+\frac{\vep}{K}\int_0^te^{-(t-\tau)/K}|u(\tau)|^{\al} d\tau,
\end{equation}
\begin{equation}\label{eq.1}
u(t)=U(t)u_0-i\int_0^t U(t-\tau)w(\tau) d\tau,
\end{equation}
with $U(t)=e^{it\De}$, $w(t)=F_0(u)(t)+F_1(u)(t)$, where
$$F_0(u)=e^{-t/K}uv_0\ \ \text{ and }\ \ 
F_1(u)=\frac{\vep}{K}u\int_0^t e^{-(t-\tau)/K}|u(\tau)|^{\al} d\tau.$$

In this setting, we show the following local well-posedness results:

\begin{maintheorem}[$n=1$]\label{t.local-A}The SDE~(\ref{e.S-D}) with cubic
nonlinearity (i.e., $\alpha=2$) is locally well-posed for $H^s\times H^s$
initial data for any $s\geq 0$. Also, the SDE~(\ref{e.S-D}) is locally
well-posed for $H^s\times H^s$ initial data when 
\begin{itemize}
\item either $s>0$ and 
$\alpha\leq 4$, or
\item $s>s_*$ and $\alpha = \frac{4}{1-2s_*}>6$.
\end{itemize} 
\end{maintheorem}

\begin{maintheorem}[$n=2$]\label{t.local-B}The SDE~(\ref{e.S-D}) with cubic
nonlinearity (i.e., $\alpha=2$) is locally well-posed for $H^s\times H^s$
initial data with $s>0$. Also, the SDE~(\ref{e.S-D}) is locally well-posed
for $H^s\times H^s$ initial data with $2\leq\alpha<\frac{2}{1-s}$.
\end{maintheorem}

\begin{maintheorem}[$n=3$]\label{t.local-C}The SDE~(\ref{e.S-D}) with cubic
nonlinearity (i.e., $\alpha=2$) is locally well-posed for $H^s\times H^s$
initial data with $s\geq 1$. Also, the SDE~(\ref{e.S-D}) is locally well-posed
for $H^s\times H^s$ initial data with $2\leq\alpha<\frac{4}{3-2s}$. 
\end{maintheorem}

\begin{maintheorem}[$n\geq 4$]\label{t.local-D}The SDE~(\ref{e.S-D}) is locally 
well-posed for $H^s\times H^s$ initial data with $2\leq\alpha<\frac{4}{n-2s}$.
\end{maintheorem}

From these local well-posedness results, the conservation of the $L^2$ norm of
$u$ and an
\emph{a priori} $H^1$ estimate we obtain the following global well-posedness
results:

\begin{maintheorem}[$n=1$]\label{t.global-E}The SDE~(\ref{e.S-D}) with cubic
nonlinearity $\alpha=2$ is globally well-posed for initial data in 
$H^s\times H^s$ with $s\geq 0$. Also, the SDE~(\ref{e.S-D}) is globally 
well-posed for $H^1\times H^1$ initial data if $\alpha\geq 1$.  
\end{maintheorem} 

\begin{maintheorem}[$n=2$]\label{t.global-F}The SDE~(\ref{e.S-D}) with cubic
nonlinearity $\alpha=2$ is globally well-posed for initial data in 
$H^s\times H^s$ with $s\geq 1$. Also, the SDE~(\ref{e.S-D}) is globally 
well-posed for $H^1\times H^1$ initial data if $\alpha\geq 2$.
\end{maintheorem}

\begin{maintheorem}[$n=3$]\label{t.global-G}The SDE~(\ref{e.S-D}) with cubic
nonlinearity $\alpha=2$ is globally well-posed for initial data in 
$H^s\times H^s$ with $s\geq 1$. Also, the SDE~(\ref{e.S-D}) is globally 
well-posed for $H^1\times H^1$ initial data if $2\leq\alpha<3$.
\end{maintheorem}

\begin{remark}A direct comparision with the global well-posedness theory of the
periodic nonlinear Schr\"odinger equation in the focusing setting~\cite{B} 
reveals a novel phenomena in the global well-posedness in the 
``focusing case'' (i.e., $\vep=-1$) of the SDE~(\ref{e.S-D}). Indeed, since the
Hamiltonian of the focusing NLS do not control the $H^1$ norm of the solutions,
we need some smallness assumptions of the $H^1$ norm of the initial data in
order to derive global wel-posedness theorems. On the other hand, the structure
of the nonlinear term of the SDE~(\ref{e.S-D}) allows us to conclude the same
global well-posedness results for the SDE \emph{without any smallness
hypothesis}. This subtle difference between the NLS and the SDE occurs because
the evolution of $v$ in the SDE permits to derive an a priori estimate for the
$H^1$ norm of $u$, although we do not have conserved Hamiltonians.
\end{remark}

We close the introduction with the scheme of this paper. In section 2, we 
revisit the restriction of the Fourier transform method of Bourgain. In 
particular, we recall the definition of the Bourgain spaces $X^{s,b}$ and some 
of its properties. Also, we revisit the Strichartz type estimates which are 
the basic tools to deal with the nonlinear terms of the SDE. In section 3, we 
prove the local well-posedness results in the theorems~\ref{t.local-A},
~\ref{t.local-B},~\ref{t.local-C} and~\ref{t.local-D}. In section 4, we derive
an a priori estimate for the $H^1$ norm. By standar arguments, this implies the
global well-posedness theorems~\ref{t.global-E},~\ref{t.global-F}
and~\ref{t.global-G}. Finally, we briefly discuss some questions related to the
well-posedness results in this paper.    

\section{Preliminaries}

This section is devoted to introduce the reader to the setting of
Bourgain~\cite{B}. 

\subsection{Restriction of the Fourier transform}

The first main ingredient of the proofs of our results is Bourgain's technique 
of restriction of the Fourier transform below. 

As Bourgain~\cite[p. 136]{B}, we are going to find a solution $u$ of the Schr\"odinger-Debye
equation~(\ref{e.S-D}) which is local in time, that is, take a function $0\leq
\psi_1\leq 1$ such that $\supp(\psi_1)\subset [-2\de,2\de]$, $\psi_1\equiv 1$ on
$[0,\de]$. In the sequel $n:=d-1$.

Inspired by equation~(\ref{eq.1}), our goal is to construct a function $u$ 
satisfying 

\begin{equation}\label{e.4-4}
u(t)=\psi_1(t)U(t)u_0-i\psi_1(t)\int_0^t U(t-\tau)w(\tau) d\tau.
\end{equation}

If we write $u_0, u, w$ as Fourier series

$$u_0(x)=\sum\limits_{\xi\in\integer^{d-1}}\hat{u_0}(\xi)
e^{2\pi i<x,\xi>}
$$

$$ u(x,t)=\sum\limits_{\xi\in\integer^{d-1}}e^{2\pi i<x,\xi>}
\int_{-\infty}^{\infty}e^{2\pi i\la t}\hat{u}(\xi,\la) \ d\la
$$

$$ w(x,t)=\sum\limits_{\xi\in\integer^{d-1}}e^{2\pi i<x,\xi>}
\int_{-\infty}^{\infty}e^{2\pi i\la t}\hat{w}(\xi,\la) \ d\la .
$$

Then the integral equation~(\ref{e.4-4}) is

\begin{eqnarray}\label{e.4-8}
u(x,t)=\psi_1(t)\sum\limits_{\xi\in\integer^{d-1}}
\hat{u_0}(\xi)e^{2\pi i(<x,\xi>+t|\xi|^2)}+\\
+\frac{1}{2\pi}\psi_1(t)\sum\limits_{\xi\in\integer^{d-1}} \nonumber
e^{2\pi i(<x,\xi>+t|\xi|^2)} \nonumber
\int_{-\infty}^{\infty}\frac{e^{2\pi i(\la-|\xi|^2)t}-1}{\la-|\xi|^2} \nonumber
\hat{w}(\xi,\la)d\la,\\ \nonumber
\end{eqnarray}

We denote by $\Phi$ the map defined by~(\ref{e.4-8}), i.e.,
\begin{eqnarray}\label{e.4-27}
\Phi(u)(t,x)=\psi_1(t)\sum\limits_{\xi\in\integer^{d-1}}
\hat{u_0}(\xi)e^{2\pi i(<x,\xi>+t|\xi|^2)}+\nonumber\\
+\frac{1}{2\pi}\psi_1(t)\sum\limits_{\xi\in\integer^{d-1}} \nonumber
e^{2\pi i(<x,\xi>+t|\xi|^2)} \nonumber
\int_{-\infty}^{\infty}\frac{e^{2\pi i(\la-|\xi|^2)t}-1}{\la-|\xi|^2} \nonumber
\hat{w}(\xi,\la)d\la,\\ \nonumber
\end{eqnarray}

Consider a function $\psi_2$ such that $0\leq\psi_2\leq 1$, $\psi_2=1$ on
$[-1,1]$ and $\supp \ \psi_2\subset [-2,2]$.

If we write,

\begin{eqnarray}
\psi_1(t)\infint\frac{e^{2\pi i(\la-|\xi|^2)t}-1}{\la-|\xi|^2} \nonumber
\hat{w}(\xi,\la)d\la=\\  \nonumber
\sum\limits_{k\geq 1}\frac{(2\pi i)^k}{k!}\psi_1(t)t^k \nonumber
\int\psi_2(\la-|\xi|^2)(\la-|\xi|^2)^{k-1}\hat{w}(\xi,\la)d\la \\ \nonumber
+\psi_1(t)\int(1-\psi_2(\la-|\xi|^2)) \nonumber
\frac{e^{2\pi i(\la-|\xi|^2)t}}{\la-|\xi|^2}\hat{w}(\xi,\la)d\la \\ \nonumber
-\psi_1(t)\int(1-\psi_2(\la-|\xi|^2)) \nonumber
\frac{\hat{w}(\xi,\la)}{\la-|\xi|^2}d\la. \nonumber
\end{eqnarray}
then the right hand side of~(\ref{e.4-8}) is controled by the 
contributions 

\begin{equation}\label{e.4-13}
\psi_1(t)\sum\limits_{\xi\in\integer^{d-1}}\hat{u_0}(\xi)e^{2\pi i(<x,\xi>+ 
t|\xi|^2)};
\end{equation}
and 
\begin{eqnarray}\label{e.4-14}
\sum\limits_{k\geq 1}\frac{(2\pi
i)^k}{k!}t^k\cdot\psi_1(t)\cdot\\ 
\left\{ \sum\limits_{\xi} 
\bigg[\int\psi_2(\la-|\xi|^2)(\la-|\xi|^2)^{k-1} \nonumber
\hat{w}(\xi,\la)\ d\la\bigg]\cdot
e^{2\pi i(<x,\xi>+ t|\xi|^2)} \right\};
\end{eqnarray}
and
\begin{equation}\label{e.4-15}
\psi_1(t)\sum\limits_{\xi\in\integer^{d-1}}e^{2\pi i<x,\xi>}\int\frac{
1-\psi_2(\la-|\xi|^2)}{\la-|\xi|^2}e^{2\pi i\la t}\hat{w}(\xi,\la)\ d\la;
\end{equation}
and
\begin{equation}\label{e.4-16}
\psi_1(t)\sum\limits_{\xi\in\integer^{d-1}}e^{2\pi i(<x,\xi>+ t|\xi|^2)}
\int\frac{1-\psi_2(\la-|\xi|^2)}{\la-|\xi|^2}\hat{w}(\xi,\la)\ d\la.
\end{equation}

We apply the Picard's fixed point method in the Bourgain spaces $X^{s,b}$ 
associated to the norm 
\begin{equation*}
\|f\|_{X^{s,b}}:=\|\langle \xi \rangle^s \langle \lambda-\xi^2 \rangle^b
\widehat{f}(\xi,\lambda)\|_{L_{\xi,\lambda}^2}.
\end{equation*}

\begin{remark}
$X^{s,b}\subset L_t^{\infty}H^s$ for any $b>1/2$. Thus, we can use the $X^{s,b}$
spaces with $b>1/2$ to prove the well-posedness results.
\end{remark}

The idea is to prove that the integral formulation of the SDE~(\ref{e.S-D}) is a
contraction of a large ball in the $X^{s,b}$. Therefore, the main task is to
estimate these four terms. Note that 
\begin{equation}\label{e.linear-homogeneous}
\|(\ref{e.4-13})\|_{X^{s,b}}\leq \|u_0\|_{H^s}.
\end{equation}
and   
\begin{equation}\label{e.linear-inhomogeneous}
\|(\ref{e.4-14})\|_{X^{s,b}}+
\|(\ref{e.4-15})\|_{X^{s,b}}+
\|(\ref{e.4-16})\|_{X^{s,b}}\leq \|w\|_{X^{s,b'-1}}
\end{equation}
for $b'>1/2$ and $b'\geq b$. Hence, it suffices to control the expression
$\|w\|_{X^{s,b'-1}}$. To accoplish this, we revisit some properties of the
Bourgain spaces.

\begin{lemma}We have 
$$\|\psi(t)f\|_{X^{s,b}}\leq_{\psi,b} 
\|f\|_{X^{s,b}}$$ 
for any $s,b\in\mathbb{R}$ and, furthermore, if $-1/2<b'\leq b <1/2$, then for 
any $0<T<1$ we have 
$$\|\psi_T(t)f\|_{X^{s,b'}}\leq_{\psi,b',b} T^{b-b'} 
\|f\|_{X^{s,b}}.$$ 
\end{lemma}  

\begin{proof}First of all, note that $\langle\lambda-\lambda_0-|\xi|^2
\rangle^{b}
\leq_b 
\langle\lambda_0\rangle^{|b|}\langle\lambda-|\xi|^2\rangle^{b}$, 
from which we obtain 
$$\|e^{it\lambda_0}f\|_{X^{s,b}}\leq_b \langle\lambda_0\rangle^{|b|} 
\|f\|_{X^{s,b}}.$$  
Using that $\psi(t)=\int\widehat{\psi}(\lambda_0) e^{it\lambda_0}d\lambda_0$, 
we conclude 
$$\|\psi(t)f\|_{X^{s,b}}\leq_b 
\left(\int|\widehat{\psi}(\lambda_0)| \langle\lambda_0\rangle^{|b|}\right) 
\|f\|_{X^{s,b}}.$$  
Since $\psi$ is smooth with compact support, the first estimate follows.  

Next we prove the second estimate. By conjugation we may assume $s=0$ and, 
by composition it suffices to treat the cases $0\leq b'\leq b$ or 
$\leq b'\leq b\leq 0$. By duality, we may take $0\leq b'\leq b$. 
Finally, by interpolation with the trivial case $b'=b$, we may consider 
$b'=0$. This reduces matters to show that 
$$\|\psi_T(t)f\|_{L^2}\leq_{\psi,b} T^b\|f\|_{X^{0,b}}$$ 
for $0<b<1/2$. Partitioning the frequency spaces into the cases 
$\langle\lambda-|\xi|^2\rangle\geq 1/T$ and $\langle\lambda-|\xi|^2\rangle
\leq 1/T$, 
we see that in the former case we'll have 
$$\|f\|_{X^{0,0}}\leq T^b\|f\|_{X^{0,b}}$$ 
and the desired estimate follows because the multiplication by $\psi$ is a 
bounded operation in Bourgain's spaces. In the latter case, by Plancherel and 
Cauchy-Schwarz 
\begin{equation*} 
\begin{split} 
\|f(t)\|_{L_x^2}&\leq \|\widehat{f(t)}(\xi)\|_{L_{\xi}^2} \leq 
\left\|\int_{\langle\lambda-|\xi|^2\rangle\leq 1/T}|\widehat{f}(\lambda,\xi)|
d\lambda) 
\right\|_{L_{\xi}^2} \\ &\leq_b T^{b-1/2} 
\left\|\int\langle\lambda-|\xi|^2\rangle^{2b} |\widehat{f}(\lambda,\xi)|^2 
d\lambda)^{1/2}\right\|_{L_{\xi}^2} = T^{b-1/2}\|f\|_{X^{s,b}}. 
\end{split} 
\end{equation*}  
Integrating this against $\psi_T$ concludes the proof of the lemma. 
\end{proof} 

In order to keep a precise control of the nonlinear term $w$, we recall the
Strichartz-type inequalities in the periodic setting derived in~\cite{B}.

\subsection{Some one-dimensional estimates}

In the 1-dimensional case, specially for the cubic nonlinearity, the 
following Strichartz estimate will be useful:

\begin{lemma}\label{l.Strichartz}
It holds $X^{0,3/8}\subset L_{x,t}^4(\mathbb{T}\times [0,1])$. More
precisely, 
\begin{equation*}
\|f\|_{L^4(\mathbb{T}\times [0,1])}\leq c\|f\|_{X^{0,3/8}}.
\end{equation*}
\end{lemma}

Next we introduce the definition:

\begin{definition}Let $d\geq 1$, $S\subset\integer^d$ and $p>2$. We define
$K_p(S)$ to be the smallest number such that
$$\left\|\sum\limits_{\ga\in S}a_{\ga} e^{2\pi i<x,\ga>}\right\|_{L^p(\Tor^d)}
\leq K_p(S)\left(\sum|a_n|^2\right)^{1/2}. $$
\end{definition}

Also, when the nonlinearity is not cubic, we will use the following
$L^6$-estimate:
\begin{proposition}\label{p.2-36}If $S_N=\{(n,n^2):|n|\leq N\}$ then
\begin{equation}\label{e.2-37}
K_6(S_N)<\exp c\frac{\log N}{\log\log N}.
\end{equation}
In particular,
\begin{equation}\label{e.2-38}
\|\sum\limits_{n\in\integer}a_ne^{i(nx+n^2t)}\|_{L^6(\Tor^2)}\ll N^{\vep}\left(
\sum|a_n|^2\right)^{1/2},\forall\ \ \vep >0.
\end{equation}
\end{proposition}

Since this proposition is not difficult to show, we include a proof of it here.

\begin{proof}Let $f=\sum\limits_{n=1}^N a_ne^{i(nx+n^2t)}$. Then:

$$\|f\|_6=\|f^3\|_2^2=\sum\limits_{n,j}\bigg|\sum\limits_{n_1^2+
n_2^2+(n-n_1-n_2)^2=j}a_{n_1}a_{n_2}a_{n-n_1-n_2}\bigg|^2. $$
Define $r_{n,j}=\#\{(n_1,n_2):|n_i|\leq N,\ n_1^2+n_2^2+(n-n_1-n_2)^2=j\}$. We
have:
$$\sum\limits_{n,j}\bigg|\sum\limits_{n_1^2+n_2^2+(n-n_1-n_2)^2=j}
a_{n_1}a_{n_2}a_{n-n_1-n_2}\bigg|^2\leq\max\limits_{|n|\leq 3N,|j|\leq 3N^2}
r_{n,j}\cdot\left(\sum\limits_{|n|\leq N}|a_n|^2\right)^3. $$

Hence, it remains to prove that
$r_{n,j}<\exp c\frac{\log N}{\log\log N}$.

The condition $n_1^2+n_2^2+(n-n_1-n_2)^2=j$ is
$n_1^2+n_2^2-nn_1-nn_2+n_1n_2=\frac{j-n^2}{2}$, i.e.,
$$\frac{3}{4}(n_1+n_2)^2+\frac{1}{4}(n_1-n_2)^2-n(n_1+n_2)=\frac{j-n^2}{2}. $$

If we put $m_1=n_1+n_2,m_2=n_1-n_2$, then:
$$(3m_1-2n)^2+3m_2^2=6j-2n^2,$$
which has the form $X^2+3Y^2=A$, where $X,Y,A\in\integer$.

Put $\rho=e^{2\pi i/3}=\frac{1+i\sqrt{3}}{2}$. Then $X^2+3Y^2=A$ if and only if 
$X+i\sqrt{3}Y$ divides $A$ in $\integer+\rho\integer$. Since 
$\integer+\rho\integer$ is an Euclidean domain, the number of divisors of $A$ is
at most $\exp c\frac{\log A}{\log\log A}< 
\exp c\frac{\log N}{\log\log N}$. 

Because $(3m_1-2n,m_2)$ defines $(n_1,n_2)$, this concludes the proof.
\end{proof}

\subsection{Some higher dimensional estimates}

For positive integers $K,N$, consider the sets
$$\La_{A,N}=\{\ze=(\xi,\la)\in\integer^n\times\real:N\leq|\xi|<2N,\ \ A\leq
|\la-|\xi|^2|<2A\}.$$

For an interval $I$ of $\integer^n$, we define
$$\La_{K,I}=\{\ze\in I\times\real:A\leq|\la-|\xi|^2|<2A\}.$$

\begin{definition}Given a function $u\in L^2(\Tor^n\times\real)$, 
$$u=\sum\limits_{\xi\in\integer^n}\int d\la \ \hat{u}(\ze) e^{2\pi 
i(<x,\xi>+\la t)}, $$
we define
\begin{equation}\label{e.5-5}
|||u|||=\sup\limits_{A,N}(A+1)^{1/2}(N+1)^s\left(\int_{\La_{A,N}}
|\hat{u}(\ze)|^2d\ze\right)^{1/2}.
\end{equation}
\end{definition}

Fix an interval $[-\de,\de]$. We will consider the restriction norm
$|||u|||=\inf|||\widetilde{u}|||$, where the infimum is taken over all
$\widetilde{u}$ coinciding with $u$ on $\Tor^n\times [0,\de]$.

Define $S_{d,N}=\{(n_1,\dots,n_{d-1},|\overline{n}|^2):n_j\in\integer,
|n_j|<N\}$, $\overline{n}=(n_1,\dots,n_{d-1})$ and
$|\overline{n}|^2=n_1^2+\dots+n_{d-1}^2$.

\begin{definition}A number $p$ is called an \emph{admissible exponent} if
\begin{equation}
p\geq\frac{2(d+1)}{d-1}\ \text{ and }\ K_p(S_{d,N})\ll
N^{\vep}N^{\frac{d-1}{2}-\frac{d+1}{p}},
\end{equation}
\end{definition}
 
Concerning the existence and the properties of 
admissible exponents we have three important results:

\begin{proposition}[Proposition 3.6 of~\cite{B}]\label{p.3-6}For $n=2,3,4$, 
the exponent $4$ is admissible, i.e.,
\begin{itemize}
\item $K_4(S_{3,N})\ll N^{\vep}$
\item $K_4(S_{4,N})\ll N^{\frac{1}{4}+\vep}$
\item $K_4(S_{5,N})\ll N^{\frac{1}{2}+\vep}$
\end{itemize}
\end{proposition}

\begin{proposition}[Proposition 3.110 of~\cite{B}]\label{p.3-110} For 
$n\geq 4$, 
$p\geq\frac{2(n+4)}{n}$,
$$K_p(S_{d,N})<cN^{\frac{d-1}{2}-\frac{d+1}{p}}. $$
\end{proposition}

\begin{proposition}[Proposition 3.113 of~\cite{B}]\label{p.3-113} If 
$p_2>p_1\geq p_0=\frac{2(d+1)}{d-1}$ and
$K_{p_1}(S_{d,N})\ll N^{\frac{d-1}{2}-\frac{d+1}{p_1}+\vep}$, then
$K_{p_2}(S_{d,N})\leq C_{p_2}N^{\frac{d-1}{2}-\frac{d+1}{p_1}}$.
\end{proposition}

The reason for the introduction of admissible exponents is explained by the good
properties (with respect to the Fourier transform) below.
  
Let $p_0$ be admissible. By proposition~\ref{p.3-113}, for $p>p_0$,
$$\left\|\sum\limits_{|\xi|\leq N}a_{\xi} e^{2\pi i(<x,\xi>+t|\xi|^2)}
\right\|_{L^p(\Tor^d)}\leq N^{\frac{d-1}{2}-\frac{d+1}{p}} 
\left(\sum\limits_{n\in\integer}|a_n|^2\right)^{1/2}.$$

Let $I$ be a $(d-1)$-interval of size $N$ in $\integer^{d-1}$ centered at
$\xi_0$. Writing
$$<x,\xi>+t|\xi|^2=<x,\xi_0>+t|\xi_0|^2+<x+2t\xi_0,\xi-\xi_0>+t|\xi-\xi_0|^2.$$

The change of variables $x'=x+2t\xi_0$, $t'=t$ implies that also
$$\left\|\sum\limits_{\xi\in I}a_{\xi} e^{2\pi i(<x,\xi>+t|\xi|^2)}
\right\|_{L^p(\Tor^d)}\leq N^{\frac{d-1}{2}-\frac{d+1}{p}} 
\left(\sum\limits_{\xi\in I}|a_{\xi}|^2\right)^{1/2}. $$
It follows that (writing $\la=|\xi|^2+k$, $|k|<A$) the map
\begin{center}
\begin{equation}\label{e.5-9}
L_{\La_{A,I}}^2\longrightarrow L^p(\Tor^n\times\real_{loc})
\end{equation}
\end{center}
\begin{center}
$\{a_{\ze}\}_{\ze\in\La_{A,I}}\to \int_{\La_{A,I}}a_{\ze} e^{2\pi
i(<x,\xi>+t\la)}d\ze$
\end{center}
has norm bounded by $A^{1/2}N^{\frac{d-1}{2}-\frac{d+1}{p}}$. 

Since the map~(\ref{e.5-9}) from $L_{\La_{A,I}}^2$ to
$L^2(\Tor^n\times\real_{loc})$ has also bounded norm, by interpolation we obtain
the following lemma:

\begin{lemma}\label{l.5-10}Let $p_1>p_0,\ p_1>p_2>2,\ \frac{1}{p_2}=
\frac{1-\theta}{p_1}+\frac{\theta}{2}$. Then the map~(\ref{e.5-9}) ranging into
$L^{p_2}(\Tor^n\times\real_{loc})$ has norm bounded by 
$$A^{\frac{1}{2}(1-\theta)}N^{(\frac{d-1}{2}-\frac{d+1}{p})(1-\theta)}.$$  
\end{lemma}

To finish these preliminaries, we introduce the notation: for a dyadic $M$,
$$u_M=\sum\limits_{|\xi|\leq M}e^{2\pi i<x,\xi>}\int\hat{u}(\xi,\la) e^{2\pi
i\la t}d\la,$$

$$\De_M u=u_M-u_{\frac{M}{2}}.$$

If $I$ is an interval of $\integer^n$,
\begin{eqnarray}
\De_I u&=&\sum\limits_{\xi\in I}e^{2\pi i<x,\xi>}\int\hat{u}(\xi,\la) 
e^{2\pi i\la t}d\la
\nonumber \\
&=&\sum\limits_{K\ dyadic}\int_{\La_{A,N}}\hat{u}(\ze) e^{2\pi i(<x,\xi>+\la t)}
d\ze.
\end{eqnarray}

This dyadic localization will be helpful in the analysis of the Bourgain norm of
the nonlinearity of the equation~(\ref{e.S-D}).

For later use, we observe that, using lemma~\ref{l.5-10}
\begin{equation}\label{e.5-35}
\|\De_I u\|_{p_2}\leq c\sum\limits_{A \ dyadic} 
A^{\frac{1}{2}(1-\theta)}M^{(\frac{d-1}{2}-\frac{d+1}{p})(1-\theta)}
\left(\int_{\La_{A,I}}|\hat{u}(\ze)|^2d\ze\right)^{1/2},
\end{equation}
if the size of $I$ is $M$.

Similarly,
\begin{eqnarray}\label{e.5-42}
\|\De_{M}u\|_{p_2}&\leq& c\sum\limits_{A \ dyadic} 
A^{\frac{1}{2}(1-\theta)}M^{(\frac{d-1}{2}-\frac{d+1}{p})(1-\theta)}
\left(\int_{\La_{A,M}}|\hat{u}(\ze)|^2d\ze\right)^{1/2}\nonumber\\
&\leq& cM^{(\frac{d-1}{2}-\frac{d+1}{p})(1-\theta)-s}\ |||u|||.
\end{eqnarray}

\section{Local well-posedness for the periodic SDE}

The basic lemma in the proof of our local well-posedness result is:

\begin{lemma}\label{l.5-64}If $s<\frac{d-1}{2}$, $\al<\frac{4}{d-1-2s}$ and 
$p_0<\frac{2(d+1)}{d-1-\frac{2}{3}s}$, where $p_0$ is an 
\emph{admissible exponent}, then 
\begin{equation}\label{e.5-65}
A^{-1/2}N^s\left(\int_{\La_{A,N}}|\hat{w}(\ze)|^2d\ze\right)^{1/2}\leq c
A^{-\theta}N^{-\theta}(\|v_0\|_{H^s}|||u|||+|||u|||^{1+\al}),
\end{equation}
for some $\theta>0$.
\end{lemma}

\begin{proof}Write $w=F_0(u)+F_1(u)$ with $F_0(u):=\mu u$ and $F_1(u):=\eta u$,
where $\mu = e^{-t/K}v_0$ and $\eta=\frac{\epsilon}{K}\int_0^t e^{-(t-\tau)/K}
|u(\tau)|^{\alpha} d\tau$. This reduces our goal to prove the estimates 
\begin{equation}\label{e.F0}
A^{-1/2}N^s\left(\int_{\La_{A,N}}|\hat{F_0}(\ze)|^2d\ze\right)^{1/2}\leq
cA^{-\theta}N^{-\theta}\|v_0\|_{H^s} |||u|||,
\end{equation}
and 
\begin{equation}\label{e.F1}
A^{-1/2}N^s\left(\int_{\La_{A,N}}|\hat{F_1}(\ze)|^2d\ze\right)^{1/2}\leq 
cA^{-\theta}N^{-\theta}|||u|||^{1+\al}.
\end{equation}

First we analyse the left-hand side of~(\ref{e.F0}). Write 
$$F_0=u\mu=e^{-t/K}\sum\left( u_M (v_0)_M - u_{\frac{M}{2}} 
(v_0)_{\frac{M}{2}}\right).$$

Hence, it suffices to prove the bound~(\ref{e.F0}) with $F_0$ replaced by 
\begin{equation*}
\Delta_M u \cdot e^{-t/K} (v_0)_M
\end{equation*}
and 
\begin{equation*}
u_{\frac{M}{2}}\cdot e^{-t/K} \Delta_M v_0 .
\end{equation*}
with $M\geq N$.

Because $u_M=\sum\limits_{M_1\leq M} \Delta_{M_1} u$,
$(v_0)_M=\sum\limits_{M_1\leq M} \Delta_{M_1} v_0$ and $\Delta_M u =
\sum\limits_{I} \Delta_I u$, where $I$ is a decomposition of $\frac{M}{2}\leq
|\xi|\leq M$ into intervals of size $M_1$, our task is to show that~(\ref{e.F0})
holds with $F_0$ replaced by 
\begin{equation*}
(F_0^{(1)})_I:=\Delta_I u\cdot e^{-t/K} \Delta_{M_1} v_0
\end{equation*}
and 
\begin{equation*}
(F_0^{(2)})_I:=\Delta_{M_1} u\cdot e^{-t/K} \Delta_{I} v_0
\end{equation*}

Choose 
$p_1>p_0$, $p_1>p_2>2$, $\frac{1}{p_2}=
\frac{1-\theta_2}{p_1}+\frac{\theta_2}{2}$.

The dual form of lemma~\ref{l.5-10} gives
\begin{equation*}
\left(\int_{\La_{A,I}}|\widehat{F_0^{(1)}}_I(\ze)|^2 d\ze\right)^{1/2}\leq 
cA^{\frac{1}{2}(1-\theta_2)}M^{(\frac{d-1}{2}-\frac{d+1}{p_1})(1-\theta_2)}
\|(F_0^{(1)})_I\|_{p_2'}
\end{equation*}
and 
\begin{equation*}
\left(\int_{\La_{A,I}}|\widehat{F_0^{(2)}}_I(\ze)|^2 d\ze\right)^{1/2}\leq 
cA^{\frac{1}{2}(1-\theta_2)}M^{(\frac{d-1}{2}-\frac{d+1}{p_1})(1-\theta_2)}
\|(F_0^{(2)})_I\|_{p_2'}
\end{equation*}
Therefore, by H\"older's inequality
\begin{equation*}
\|(F_0^{(1)})_I\|_{p_2'}\leq \frac{1}{K}\|\De_I u\|_{p_2}
\| e^{-t/K}\De_{M_1}v_0\|_{\frac{p_2-p_2'}{p_2p_2'}}.
\end{equation*}
and 
\begin{equation*}
\|(F_0^{(2)})_I\|_{p_2'}\leq \frac{1}{K}\|\De_I v_0\|_{p_2}
\| e^{-t/K}\De_{M_1}u\|_{\frac{p_2-p_2'}{p_2p_2'}}.
\end{equation*}
Using~(\ref{e.5-35}), we obtain
\begin{equation*}
\left(\sum\limits_I\|\De_I u\|_{p_2}^2\right)\leq
cM_1^{(\frac{d-1}{2}-\frac{d+1}{p_1})(1-\theta_2)}\cdot M^{-s}|||u|||.
\end{equation*}
and 
\begin{equation*}
\left(\sum\limits_I\|\De_I v_0\|_{p_2}^2\right)\leq
cM_1^{(\frac{d-1}{2}-\frac{d+1}{p_1})(1-\theta_2)}\cdot M^{-s}\|v_0\|_{H^s}.
\end{equation*}
Also, taking 
\begin{equation*}
p_3>p_0, p_3>p_4>2,
\frac{1}{p_4}=\frac{1-\theta_4}{p_3}+\frac{\theta_4}{2} \text{ and } 1> 
\frac{2}{p_2}+
\frac{1}{p_4}
\end{equation*}
then, the estimate~(\ref{e.5-42}) implies
\begin{equation*}
\| e^{-t/K}\De_{M_1}u\|_{\frac{p_2-p_2'}{p_2p_2'}}\leq 
\|\De_{M_1}u\|_{p_4}\leq M_1^{(\frac{d-1}{2}-\frac{d+1}{p_3})(1-\theta_4)-s}
|||u|||
\end{equation*}
and 
\begin{equation*}
\| e^{-t/K}\De_{M_1}v_0\|_{\frac{p_2-p_2'}{p_2p_2'}}\leq 
\|\De_{M_1}v_0\|_{p_4}\leq M_1^{(\frac{d-1}{2}-\frac{d+1}{p_3})(1-\theta_4)-s}
\|v_0\|_{H^s}.
\end{equation*}
Thus, after performing the summations over $M_1\leq M$ and $M\geq N$ in the
previous estimates, we get the desired bounds on $F_0^{(1)}$ and $F_0^{(2)}$. In
particular,~(\ref{e.F0}) is proved, if there are numbers $p_1,\dots, p_4$
verifying the relations above.

Next, we show the estimate~(\ref{e.F1}). Write
$$F_1=u\eta=\frac{\vep}{K}\sum\left( u_M\int_0^t e^{-(t-\tau)/K}
|u_M|^{\al}d\tau-u_{\frac{M}{2}}\int_0^t e^{-(t-\tau)/K}
|u_{\frac{M}{2}}|^{\al}d\tau\right).$$

So, we have to
evaluate~(\ref{e.F1}) with $F_1$ replaced by 
$$\frac{\vep}{K}\De_M u\cdot\int_0^t
e^{-(t-\tau)/K}|u_M|^{\al}d\tau,
$$
and
$$
\frac{\vep}{K}u_{\frac{M}{2}}
\int_0^t e^{-(t-\tau)/K}(|u_M|^{\al}-|u_{\frac{M}{2}}|^{\al})d\tau.
$$
with $M\geq N$.

Since for $\al\geq 2$ and complex numbers $z,w$,
$$|z|^{\al}-|w|^{\al}=(z-w)\phi_1(z,w)+(\overline{z}-\overline{w})\phi_2(z,w),$$
where $|\phi_1|,|\phi_2|\leq c(|z|+|w|)^{\al-1}$, if we write
$u_M=\sum\limits_{M_1\leq M}\De_{M_1}u$, it is sufficient to
estimate~(\ref{e.F1}) with $F_1$ replaced by 
\begin{equation}\label{e.F1a}
\frac{\vep}{K}\De_M u\cdot\int_0^t
e^{-(t-\tau)/K}\De_{M_1}u\cdot\phi(u_{M_1},u_{\frac{M_1}{2}})d\tau
\end{equation} 
and
\begin{equation}\label{e.F1b}
\frac{\vep}{K}\De_{M_1} u\cdot\int_0^t
e^{-(t-\tau)/K}\De_{M}u\cdot\phi(u_{M},u_{\frac{M}{2}})d\tau
\end{equation} 
where $M_1\leq M$, $M\geq N$. We subdivide $\frac{M}{2}<|\xi|\leq M$ in 
intervals
$I$ of size $M_1$ and write 
$$\De_M u=\sum\limits_I\De_I u.$$

Because the functions $\SA_I=\frac{\vep}{K}\De_I u\cdot\int_0^t
e^{-(t-\tau)/K}\De_{M_1}u\cdot\phi(u_{M_1},u_{\frac{M_1}{2}})d\tau$ (resp.
$\SB_I=\frac{\vep}{K}\De_{M_1} u\cdot\int_0^t
e^{-(t-\tau)/K}\De_{I}u\cdot\phi(u_{M},u_{\frac{M}{2}})d\tau$) have
essentially disjoint supports, the contributions of~(\ref{e.F1a}),~(\ref{e.F1b})
to~(\ref{e.F1}) are

\begin{equation}\label{e.5-30a}
A^{-1/2}N^s\left(\sum\limits_I\int_{\La_{A,I}}|\hat{\SA}_I(\ze)|^2 
d\ze\right)^{1/2}
\end{equation}
and
\begin{equation}\label{e.5-30b}
A^{-1/2}N^s\left(\sum\limits_I\int_{\La_{A,I}}|\hat{\SB}_I(\ze)|^2 
d\ze\right)^{1/2}.
\end{equation}

We deal first with the contribution~(\ref{e.5-30a}). Choose 
$p_1>p_0$, $p_1>p_2>2$, $\frac{1}{p_2}=
\frac{1-\theta_2}{p_1}+\frac{\theta_2}{2}$.

The dual form of lemma~\ref{l.5-10} gives
\begin{equation}\label{e.5-32}
\left(\int_{\La_{A,I}}|\hat{\SA}_I(\ze)|^2 d\ze\right)^{1/2}\leq 
cA^{\frac{1}{2}(1-\theta_2)}M^{(\frac{d-1}{2}-\frac{d+1}{p_1})(1-\theta_2)}
\|\SA_I\|_{p_2'}
\end{equation}
and, by H\"older's inequality
\begin{equation}\label{e.5-33}
\|\SA_I\|_{p_2'}\leq \frac{1}{K}\|\De_I u\|_{p_2}
\|\int_0^t e^{-(t-\tau)/K}\De_{M_1}u\cdot\phi(u_{M_1},u_{\frac{M_1}{2}})
d\tau\|_{\frac{p_2-p_2'}{p_2p_2'}}.
\end{equation}

Using~(\ref{e.5-35}), we obtain
\begin{equation}\label{e.5-37}
\left(\sum\limits_I\|\De_I u\|_{p_2}^2\right)\leq
cM_1^{(\frac{d-1}{2}-\frac{d+1}{p_1})(1-\theta_2)}\cdot M^{-s}|||u|||.
\end{equation}

On the other hand, choosing 
\begin{equation}\label{e.5-39}
p_3>p_0, p_3>p_4>2,
\frac{1}{p_4}=\frac{1-\theta_4}{p_3}+\frac{\theta_4}{2} \text{ and } 1> 
\frac{2}{p_2}+
\frac{1}{p_4}
\end{equation}
then
\begin{equation}\label{e.5-40}
\|\int_0^t e^{-(t-\tau)/K}\De_{M_1}u\cdot\phi(u_{M_1},u_{\frac{M_1}{2}})
d\tau\|_{\frac{p_2-p_2'}{p_2p_2'}}\leq \|\De_{M_1}u\|_{p_4}\cdot 
\|\phi\|_{(1-\frac{2}{p_2}-
\frac{1}{p_4})^{-1}}
\end{equation}

Note that
\begin{equation}\label{e.5-43}
\|\phi\|_{(1-\frac{2}{p_2}-\frac{1}{p_4})^{-1}}\leq
\|u_{M_1}\|_{(\al-1)(1-\frac{2}{p_2}-\frac{1}{p_4})^{-1}}^{\al-1}
\end{equation}

But, writing $u_{M_1}=\sum\limits_{M_2< M_1 \ dyadic}\De_{M_2}u$, if we choose
$p_5>p_0, p_5>p_6>2$, $\frac{1}{p_6}=\frac{1-\theta_6}{p_5}+\frac{\theta_6}{2}$ 
and $\frac{\al-1}{p_6}\leq 1-\frac{2}{p_2}-\frac{1}{p_4}$, then 

\begin{equation}\label{e.5-45}
(\ref{e.5-43})\leq c|||u|||^{\al-1}.
\end{equation}

Putting together the estimates~(\ref{e.5-45}),~(\ref{e.5-40}),~(\ref{e.5-42}) 
and
performing summations over $M_1\leq M$ and $M\geq N$, we proved
\begin{equation}\label{e.5-50a}
(\ref{e.5-30a})\leq cA^{-\theta}N^{-\theta}\cdot |||u|||^{1+\al}
\end{equation}
for some $\theta>0$, provided that we can assure the existence of $p_1,\dots,
p_6$ satisfying the relations above. 

Similarly, the contribution of~(\ref{e.5-30b}) can be analysed as follows.
Keeping the same notation as above, the dual form of the lemma~\ref{l.5-10}
still yields 
\begin{equation}\label{e.5-32b}
\left(\int_{\La_{A,I}}|\hat{\SB}_I(\ze)|^2 d\ze\right)^{1/2}\leq 
cA^{\frac{1}{2}(1-\theta_2)}M^{(\frac{d-1}{2}-\frac{d+1}{p_1})(1-\theta_2)}
\|\SB_I\|_{p_2'}
\end{equation}
and, by H\"older's inequality
\begin{equation}\label{e.5-33b}
\|\SB_I\|_{p_2'}\leq \frac{1}{K}\|\De_{M_1} u\|_{p_4}
\|\int_0^t e^{-(t-\tau)/K}\De_{I}u\cdot\phi(u_{M_1},u_{\frac{M_1}{2}})
d\tau\|_{\widetilde{p}_4},
\end{equation}
where $\widetilde{p}_4=(1-\frac{1}{p_2}-\frac{1}{p_4})^{-1}$.

Using~(\ref{e.5-35}), we obtain again 
\begin{equation}\label{e.5-37b}
\left(\sum\limits_I\|\De_I u\|_{p_2}^2\right)\leq
cM_1^{(\frac{d-1}{2}-\frac{d+1}{p_1})(1-\theta_2)}\cdot M^{-s}|||u|||.
\end{equation}

On the other hand, choosing 
\begin{equation}\label{e.5-39b}
\frac{1}{\widetilde{p}_4}=\frac{1}{p_2}+\frac{1}{\widehat{p}_4}
\end{equation}
then, since $\widehat{p}_4=(1-\frac{2}{p_2}-
\frac{1}{p_4})^{-1}$, 
\begin{equation}\label{e.5-40b}
\|\int_0^t e^{-(t-\tau)/K}\De_{M_1}u\cdot\phi(u_{M_1},u_{\frac{M_1}{2}})
d\tau\|_{\widetilde{p}_4}\leq \|\De_{I}u\|_{p_2}\cdot 
\|\phi\|_{(1-\frac{2}{p_2}-
\frac{1}{p_4})^{-1}}
\end{equation}
Thus, we can apply the same arguments used in the treatment of~(\ref{e.5-30a})
to get 
\begin{equation}\label{e.5-50b}
(\ref{e.5-30b})\leq c A^{-\theta}N^{-\theta}|||u|||^{1+\alpha}.
\end{equation}

Finally, it remains only to justify the existence of the numbers $p_1,\dots,
p_6$ satisfying the claimed relations. However, it is not difficult to prove 
(see~\cite[p.149]{B}) that these numbers exist if $s<\frac{d-1}{2}$, 
$\al<\frac{4}{d-1-2s}$ and $p_0<
\frac{2(d+1)}{d-1-\frac{2}{3}s}$.
\end{proof}

Once this lemma is proved, it is a standard matter to get the local
well-posedness statements in the theorems~\ref{t.local-A},~\ref{t.local-B},
~\ref{t.local-C} and~\ref{t.local-D}. Indeed, the lemma~\ref{l.5-64} can be
applied to give the estimate 
\begin{equation*}
\|w\|_{X^{s,-1/2+}}\leq c (\|v_0\|_{H^s}\|u\|_{X^{s,1/2}}+
\|u\|_{X^{s,1/2}}^{1+\alpha}).
\end{equation*}
In particular, this estimate can be combined with the
bounds~(\ref{e.linear-homogeneous}) and~(\ref{e.linear-inhomogeneous}) to obtain
that the integral formulation of the SDE~(\ref{e.S-D}) is a contraction of a 
large ball in the space $X^{s,b}$ into itself. This completes the proof of the
local well-posedness thereoms~\ref{t.local-A},~\ref{t.local-B},
~\ref{t.local-C} and~\ref{t.local-D}. 

\section{Global well-posedness for the periodic SDE}

We start with the case of cubic nonlinearity in dimensions $n=1, 2, 3$: the 
proof of the theorem~\ref{t.global-E} clearly follows from the conservation of the $L^2$-norm
of $u$, if we can prove the estimate 
\begin{equation}\label{e.cubic-global-1}
\|w\|_{X^{s,-1/2+}}\leq c(\|v_0\|_{H^s}\|u\|_{X^{s,1/2}}+ \|u\|_{X^{0,1/2}}^2
\|u\|_{X^{s,1/2}}).
\end{equation}  
Similarly, the proof of theorems~\ref{t.global-F} and~\ref{t.global-G} follows
from the estimate 
\begin{equation}\label{e.cubic-global-2}
\|w\|_{X^{s,-1/2+}}\leq c(\|v_0\|_{H^s}\|u\|_{X^{s,1/2}}+\|u\|_{X^{1,1/2}}^2
\|u\|_{X^{s,1/2}}).
\end{equation}
However, the bound in~(\ref{e.cubic-global-1}) is easily obtained via a simple
modification of the calculations in~\cite[p.110--114]{B3} using the Strichartz
estimate in lemma~\ref{l.Strichartz}. Analogously, the
bound~(\ref{e.cubic-global-2}) follows from simple modifications of the
calculations in~\cite[p.115--118]{B3} (along the lines of the proof of the
lemma~\ref{l.5-64}) using the Strichartz bounds in propositions~\ref{p.3-6}
and~\ref{p.3-110}.
 
Next, we study the 
variation of the $H^1$-norm of $u$ (see the proposition
below). Using this, we will derive an \emph{a priori estimate} for the solution.
\begin{proposition}\label{p.conserv}
\begin{equation}\frac{d}{dt}\left(\int_{\Tor^n}|\nabla u(t)|^2-
\int_{\Tor^n}|u(t)|^2 v(t)\right)= \frac{1}{K}\cdot\left(\int_{\Tor^n}|u(t)|^2
v(t)-\vep\int_{\Tor^n}|u(t)|^p\right).
\end{equation}
\end{proposition}

\begin{proof} Write $u=a+i b$. The equation~(\ref{e.S-D}) implies that
\begin{eqnarray}\label{e.S-D1}
\left\{ \begin{array}{ll}
\p_t a=-\De b+bv,\\
\p_t b=\De a-av
\end{array}
\right.
\end{eqnarray}
But,
\begin{eqnarray}
\frac{1}{2}\frac{d}{dt}\int_{\Tor^n}|\nabla u(t)|^2&=&
\int_{\Tor^n}<\nabla a,\nabla\p_t a>+\int_{\Tor^n}<\nabla b,\nabla\p_t b>=
\nonumber\\
&=&\int_{\Tor^n}\left(\p_t a\De a+\p_t b\De b\right)\nonumber.
\end{eqnarray}
Hence by equation~(\ref{e.S-D1}),
\begin{equation}\label{e.conserv1}
\frac{1}{2}\frac{d}{dt}\int_{\Tor^n}|\nabla u(t)|^2=
\int_{\Tor^n}\left(b\De av-a\De bv\right).
\end{equation}
On the other hand, the equation~(\ref{e.S-D}) also implies
\begin{equation}\label{e.S-D2}
\p_t v=-\frac{1}{K}v+\frac{\vep}{K}|u|^{\al}
\end{equation}
However,
\begin{eqnarray}
\frac{1}{2}\frac{d}{dt}\int_{\Tor^n}|u(t)|^2 v(t)&=&
\frac{1}{2}\int_{\Tor^n}v(t)\p_t|u|^2+ \frac{1}{2}\int_{\Tor^n}|u(t)|^2\p_t
v\nonumber.
\end{eqnarray}
So using equations~(\ref{e.S-D1}),~(\ref{e.S-D2}), we have
\begin{equation}\label{e.conserv2}
\frac{1}{2}\frac{d}{dt}\int_{\Tor^n}|u(t)|^2 v(t)=
\int_{\Tor^n}\left(b\De av-a\De bv\right)
-\frac{1}{2K}\int_{\Tor^n}|u(t)|^2 v(t)+\frac{\vep}{2K}\int_{\Tor^n}|u(t)|^p.
\end{equation}
Then, if we subtract the equations~(\ref{e.conserv1}) and~(\ref{e.conserv2}),
the proof is complete.
\end{proof}

Integrating the equation of proposition~\ref{p.conserv}, we obtain
\begin{equation}\label{e.6-5}
\int_{\Tor^n}|\nabla u(t)|^2=\int_{\Tor^n}|u(t)|^2v(t) 
-\int_{\Tor^n}|\nabla u_0|^2+\int_{\Tor^n}|u_0|^2v_0 
+\frac{1}{K}\int|u|^2v-\frac{\vep}{K}\int|u|^p.
\end{equation}

We recall the following basic inequality
\begin{equation}\label{e.interp}
\|f\|_{L^p(\Tor^n)}\leq c\|f\|_2^{1-\theta}\|f\|_{H^1}^{\theta},
\end{equation}
where $\theta:=n\left(\frac{1}{2}-\frac{1}{p}\right)<1$.

Then, by H\"older inequality,
\begin{eqnarray}
\int_{\Tor^n}|u(t)|^2v(t)\leq\|u(t)\|_4^2\|v(t)\|_2,\nonumber
\end{eqnarray}
\begin{eqnarray}
\int|u|^2v\leq \|u\|_4^2\|v\|_2\nonumber
\end{eqnarray}
But, by~(\ref{e.interp}), since $\|u(t)\|_2=\|u_0\|_2$,
\begin{eqnarray}
\|u(t)\|_4\leq c\|u_0\|_2^{1-\theta_0}\|u(t)\|_{H^1}^{\theta_0}\nonumber,
\end{eqnarray}
\begin{eqnarray}
\|u\|_4\leq cT^{1/4}\|u_0\|_2^{1-\theta_0}
\sup\limits_{t\in [0,T]}\|u(t)\|_{H^1}^{\theta_0},\nonumber
\end{eqnarray}
\begin{eqnarray}
\int|u|^p\leq cT\|u_0\|_2^{p(1-\theta)}
\sup\limits_{t\in [0,T]}\|u(t)\|_{H^1}^{p\theta}.\nonumber
\end{eqnarray}
with $\theta_0=n(\frac{1}{2}-\frac{1}{4})=\frac{n}{4}<1$,
$\theta=n(\frac{1}{2}-\frac{1}{p})<1$. Moreover,
\begin{eqnarray}
\|v(t)\|_2\leq\|v_0\|_2+\frac{c}{K}\|u\|_{2\al}^{\al}\leq\nonumber\\
\|v_0\|_2+\frac{c}{K}T^{1/2}\|u_0\|_2^{\al(1-\theta_1)}
\sup\limits_{t\in [0,T]}\|u(t)\|_{H^1}^{\al\theta_1}\nonumber,
\end{eqnarray}
\begin{eqnarray}
\|v\|_2\leq T^{1/2}\|v_0\|_2+\frac{c}{K}T\|u_0\|_2^{\al(1-\theta_1)}
\sup\limits_{t\in [0,T]}\|u(t)\|_{H^1}^{\al\theta_1}\nonumber.
\end{eqnarray}
where $\theta_1=n(\frac{1}{2}-\frac{1}{2\al})<1$.

Applying these inequalities for equation~(\ref{e.6-5}), we get the following
\emph{a priori} estimate:

\begin{eqnarray}\label{e.apriori}
\sup\limits_{t\in [0,T]}\|u(t)\|_{H^1}^2\leq\|u_0\|_{H^1}^2+
c\|v_0\|_2\|u_0\|_2^{\frac{(4-n)}{2}}\sup\limits_{t\in
[0,T]}\|u(t)\|_{H^1}^{\frac{n}{2}}+\nonumber\\
c\|u_0\|_2^{\frac{(4-n)}{2}}\mu_1(T)\sup\limits_{t\in
[0,T]}\|u(t)\|_{H^1}^{\frac{n}{2}+\al\theta_1}+\nonumber\\
\frac{c}{K}T\|u_0\|_2^{\frac{(4-n)}{2}}\left(\|v_0\|+\mu_1(T)\sup\limits_{t\in
[0,T]}\|u(t)\|_{H^1}^{\al\theta_1}\right)\sup\limits_{t\in
[0,T]}\|u(t)\|_{H^1}^{\frac{n}{2}}+\nonumber\\
\frac{c}{K}T\|u_0\|_2^{p(1-\theta)}
\sup\limits_{t\in [0,T]}\|u(t)\|_{H^1}^{p\theta}.
\end{eqnarray}
where $\mu_1(T)=\frac{c}{K}T^{1/2}\|u_0\|_2^{\al(1-\theta_1)}$. From the
previous a priori estimate, using a standard argument, 
if $\theta_0,\theta_1,\theta<1$, then we will obtain our global well-posedness
results in the theorems~\ref{t.global-E},~\ref{t.global-F} and~\ref{t.global-G} 
for $H^1\times H^1$ data, as follows. 

Note that $\theta_0<1\iff n\leq 3$. Also, if $n=1,2$, $\theta_1<1$
for any $\al>0$ (i.e., any $p$); if $n=3$, $\theta_1<1\iff \al<3$ 
(i.e., $p<5$). Finally, $\theta<1\iff p<\frac{2n}{n-2}$. These informations 
together clearly gives the desired results.

\section{Concluding remarks}

We finish this article with two questions motivated by the previous results.
Firstly, in view of the global well-posedness theorem for the periodic NLS
equation in dimension $4$ proved by Bourgain in~\cite{B2}, it is natural to ask:

\begin{question}In dimension $4$, is the periodic SDE~(\ref{e.S-D}) where the 
nonlinearity $|u|^{\alpha}$ is replaced by $f(|u|^2)$ with 
$f(t)=O''(t^{1/2})$ (i.e., $|f(t)|\leq c t^{1/2}$, $|f'(t)|\leq c t^{-1/2}$ and 
$|f''(t)|\leq c t^{-3/2}$) globally well-posed for $H^s\times H^s$ initial data
satisfying $s\geq 2$?
\end{question}

Secondly, while our results are always stated for $H^s\times H^s$ initial data, 
Corcho and Linares~\cite{CL} were able to prove well-posedness for $H^k\times
H^s$ initial data with $k\neq s$. Thus, a interesting question is: 

\begin{question}Is the periodic SDE~(\ref{e.S-D}) well-posed for $H^k\times H^s$
initial data with $k\neq s$?
\end{question}

We plan to attack these issues in forthcoming papers. At the present moment, we 
advance that some work in progress by Corcho and the second indicates the 
possibitity of a satisfactory answer for the second question in dimension $1$.

\bibliographystyle{alpha}
\bibliography{bib}


\vspace{2.5cm}

\noindent    \textbf{Alexander Arbieto} ( alexande{\@@}impa.br )\\
	     \textbf{Carlos Matheus} ( matheus{\@@}impa.br )\\
             IMPA, Est. D. Castorina 110, Jardim Bot\^anico, 22460-320 \\
             Rio de Janeiro, RJ, Brazil

\end{document}